\input amstex
\input amsppt.sty \magnification=1200
\NoBlackBoxes \hsize=16.8truecm \vsize=22.6truecm
\def\q{\quad}
\def\qq{\qquad}
\def\qtq#1{\q\t{#1}\q}

\par\q
\hbox{}\par\q\par
\def\mod#1{\ (\text{\rm mod}\ #1)}
\def\t{\text}
\def\f{\frac}

\def\e{\equiv}
\def\b{\binom}
\def\ap{\langle a\rangle_p}

\def\sls#1#2{(\f{#1}{#2})}

\def\Ls#1#2{\Big(\f{#1}{#2}\Big)}
\let \pro=\proclaim
\let \endpro=\endproclaim

\topmatter
\title {Super congruences involving Bernoulli and Euler polynomials}\endtitle
\author Zhi-Hong Sun \endauthor
\affil School of Mathematical Sciences, Huaiyin Normal University,
\\ Huaian, Jiangsu 223001, PR China
\\ Email: zhihongsun$\@$yahoo.com
\\ Homepage: http://www.hytc.edu.cn/xsjl/szh
\endaffil

 \nologo \NoRunningHeads
 \nologo \NoRunningHeads

\abstract{Let $p>3$ be a prime, and let $a$ be a rational p-adic
integer.  Let $\{B_n(x)\}$ and $\{E_n(x)\}$ denote the Bernoulli
polynomials and Euler polynomials, respectively. In this paper we
show that
$$\sum_{k=0}^{p-1}\binom ak\binom{-1-a}k\equiv
(-1)^{\langle a\rangle_p}+ p^2t(t+1)E_{p-3}(-a)\pmod{p^3}$$ and for
$a\not\equiv -\frac 12\pmod p$,
$$\sum_{k=0}^{p-1}\binom ak\binom{-1-a}k\frac 1{2k+1}\equiv
\frac{1+2t}{1+2a} +p^2\frac{t(t+1)}{1+2a}B_{p-2}(-a)\pmod{p^3},$$
where $\langle a\rangle_p\in\{0,1,\ldots,p-1\}$ satisfying $a\equiv
\langle a\rangle_p\pmod p$ and $t=(a-\langle a\rangle_p)/p$. Taking
$a=-\frac 13,-\frac 14,-\frac 16$ in the above congruences we solve
some conjectures of Z.W. Sun. In this paper we also establish
congruences for $\sum_{k=0}^{p-1}k\binom ak\binom{-1-a}k,\
\sum_{k=0}^{p-1}\binom ak\binom{-1-a}k\frac 1{2k-1},\
\sum_{k=1}^{p-1}\frac 1k\binom ak\binom{-1-a}k\pmod{p^3}$ and
$\sum_{k=1}^{p-1}\frac {(-1)^k}k\binom ak,\ \sum_{k=0}^{p-1}\binom
ak(-2)^k\pmod{p^2}.$

\par\q
}
 \endabstract
 \footnote"" {The author is
supported by the National Natural Science Foundation of China (grant
no. 11371163).}
 \endtopmatter

\document
\subheading{1. Introduction}
\par Let $p>3$ be a prime. In 2003, based on his work concerning
hypergeometric functions and Calabi-Yau manifolds,
Rodriguez-Villegas [RV] conjectured the following congruences:
$$\align &\sum_{k=0}^{p-1}\f{\b{2k}k^2}{16^k}\e
\Ls{-1}p\mod{p^2},\tag 1.1
\\&\sum_{k=0}^{p-1}
\f{\b{2k}k\b{3k}k}{27^k}\e\Ls {-3}p\mod{p^2},\tag 1.2
\\&\sum_{k=0}^{p-1}\f{\b{2k}k\b{4k}{2k}}{64^k}\e \Ls{-2}p\mod{p^2},
\tag 1.3
\\& \sum_{k=0}^{p-1} \f{\b{3k}k\b{6k}{3k}}{432^k}\e
\Ls{-1}p\mod{p^2},\tag 1.4\endalign$$ where $\sls ap$ is the
Legendre symbol. These congruences were later confirmed by Mortenson
[M1-M2] via the Gross-Koblitz formula. For elementary proofs of
(1.1) see [S5] and [T1]. For elementary proofs of (1.2)-(1.4) see
[S8].
\par The Bernoulli numbers $\{B_n\}$ and
Bernoulli polynomials $\{B_n(x)\}$ are
 defined by
 $$B_0=1,\ \sum_{k=0}^{n-1}\b nkB_k=0\ (n\ge 2)\qtq{and}
 B_n(x)=\sum_{k=0}^n\b nkB_kx^{n-k}\ (n\ge 0).$$
 The Euler numbers $\{E_n\}$ and
Euler polynomials $\{E_n(x)\}$ are
 defined by
$$\align &E_0=1,\ E_n=-\sum_{k=1}^{[n/2]}\b n{2k}E_{n-2k}\ (n\ge 1)
\\&\t{and}\q E_n(x)=\f 1{2^n}\sum_{k=0}^n\b nk(2x-1)^{n-k}E_k,
\endalign$$
 where $[a]$ is the greatest integer not exceeding $a$.
 It is well
 known that $B_{2n+1}=0$ and $E_{2n-1}=0$ for any positive integer
 $n$. $\{B_n\}$ and $\{E_n\}$ are important sequences and they have many
interesting properties and applications. See [B], [MOS] and
[S1,S2,S3,S4].
\par Let $p>3$ be a prime and $H_n=1+\f 12+\cdots+\f 1n$.  In [Su1], using a complicated method Z.W. Sun proved that
$$\sum_{k=0}^{p-1}\f{\b{2k}k^2}{16^k}\e
\Ls{-1}p-p^2E_{p-3}\mod{p^3}\tag 1.5$$ and conjectured that (see
[Su1, Conjecture 5.12] and [Su2, Conjecture 1.2])
$$\align &\sum_{k=0}^{p-1}\f{\b{6k}{3k}\b{3k}k}{432^k}
\e \Ls {-1}p-\f {25}9p^2E_{p-3}\mod{p^3},\tag 1.6
\\&\sum_{k=0}^{p-1}\f{\b{2k}{k}\b{4k}{2k}}{64^k}
\e \Ls {-2}p-\f 3{16}p^2E_{p-3}\Big(\f 14\Big)\mod{p^3},\tag 1.7
\\&\sum_{k=0}^{p-1}\f{\b{2k}k\b{3k}k}{27^k}
\e \Ls {-3}p-\f {p^2}3B_{p-2}\Ls 13\mod{p^3},\tag 1.8
\\&\sum_{k=0}^{p-1}\f{\b{2k}k\b{4k}{2k}}{64^k(2k+1)}
\e \Ls{-1}p-3p^2E_{p-3}\mod{p^3},\tag 1.9
\\&\sum_{k=0}^{p-1}\f{\b{2k}k\b{4k}{2k}}{64^kk}\e -3H_{\f{p-1}2}+
\f 74p^2B_{p-3}\mod{p^3}.\tag 1.10
\endalign$$
\par As pointed out in [S8], we have
$$\aligned&\b{-\f 12}k^2=\f{\b{2k}k^2}{16^k},\ \b{-\f
13}k\b{-\f 23}k=\f{\b{2k}k\b{3k}k}{27^k}, \\&\ \b{-\f 14}k\b{-\f
34}k=\f{\b{2k}k\b{4k}{2k}}{64^k},\ \b{-\f 16}k\b{-\f
56}k=\f{\b{3k}k\b{6k}{3k}}{432^k}.\endaligned\tag 1.11$$
 Let $\Bbb
Z$ be the set of integers. For a prime $p$ let $\Bbb Z_p$ denote the
set of rational $p-adic$ integers. For a $p-adic$ integer $a$ let
$\ap\in\{0,1,\ldots,p-1\}$ be given by $a\e\ap\mod p$. Let $p$ be an
odd prime and $a\in\Bbb Z_p$. In [S8] the author showed that
$$\sum_{k=0}^{p-1}\b ak\b{-1-a}k\e (-1)^{\ap}\mod{p^2}.
\tag 1.12$$ For $a=-\f 12,-\f 13,-\f 14,-\f 16$, using (1.11) we get
(1.1)-(1.4) immediately. In [T3] Tauraso obtained a congruence for
$\sum_{k=1}^{p-1}\f 1k\b ak\b{-1-a}k\mod{p^2}$.
\par For a prime $p>3$ and $a\in\Bbb Z_p$ with $a\not\e 0\mod p$, in
Section 2 we improve (1.12) by showing that
$$\sum_{k=0}^{p-1}\b ak\b{-1-a}k\e (-1)^{\ap}
+(a-\ap)(p+a-\ap)E_{p-3}(-a)\mod{p^3}.\tag 1.13$$ Taking $a=-\f
12,-\f 13,-\f 14,-\f 16$ in (1.13) we deduce (1.5)-(1.8).
\par Let $p>3$ be a prime and $a\in\Bbb Z_p$ with $a\not\e -\f
12\mod p$. In Section 3 we prove that
$$\sum_{k=0}^{p-1}\b ak\b{-1-a}k\f 1{2k+1}\e \f{1+2t}{1+2a}
+p^2\f{t(t+1)}{1+2a}B_{p-2}(-a)\mod{p^3},\tag 1.14$$ where
$t=(a-\ap)/p$. Taking $a=-\f 14$ in (1.14) we deduce (1.9). In
Section 4 we determine $\sum_{k=0}^{p-1}\b ak\b{-1-a}k\f 1{2k-1}
\mod{p^3}$. In Section 5 we give a congruence for
$\sum_{k=1}^{p-1}\f 1k\b ak\b{-1-a}k$ $\mod {p^3}$. By taking $a=-\f
14$ we get (1.10). In Section 6 we give a congruence for
$\sum_{k=1}^{p-1}\f {(-1)^k}k\b ak\mod {p^2}$, and in Section 7 we
establish a congruence for $\sum_{k=0}^{p-1}\b ak(-2)^k$ $\mod
{p^2}$.

\subheading{2. Congruences for $\sum_{k=0}^{p-1}\b ak\b{-1-a}k\mod
{p^3}$}
 \pro{Lemma 2.1} Let $p>3$ be a prime
and $t\in\Bbb Z_p$. Then
$$\sum_{k=0}^{p-1}\b{pt}k\b{-1-pt}k
\e 1\mod{p^3}.$$
\endpro
Proof. For $k\in\{1,2,\ldots,p-1\}$ we see that
$$\align \b{pt}k\b{-1-pt}k&=\f{pt(pt-1)\cdots(pt-k+1)
(-1-pt)(-2-pt)\cdots(-k-pt)}{k!^2}
\\&=\f{(-1)^kpt(pt+k)}{k!^2}(p^2t^2-1^2)\cdots
(p^2t^2-(k-1)^2)
\\&\e -\f{pt(pt+k)}{k^2}=-\f{p^2t^2}{k^2}-\f{pt}k\mod{p^3}.\endalign$$
From [L] or [S2] we know that
$$\sum_{k=1}^{p-1}\f 1{k^2}\e 0\mod p\qtq{and}\sum_{k=1}^{p-1}\f 1k
\e 0\mod{p^2}.\tag 2.1$$
 Thus,
$$\sum_{k=0}^{p-1}\b{pt}k\b{-1-pt}k
\e 1-p^2t^2\sum_{k=1}^{p-1}\f 1{k^2}-pt\sum_{k=1}^{p-1}\f 1k \e
1\mod{p^3}.$$ This proves the lemma.

 \pro{Lemma 2.2} Let $p$ be an odd prime, $a\in\Bbb Z_p$,
  $a\not\e 0\mod p$ and $k\in\{1,2,\ldots,p-2\}$. Then
  $$\sum_{r=1}^{\ap}\f{(-1)^r}{r^k}
  \e -\f{(2^{p-k}-1)B_{p-k}}{p-k}+\f 12(-1)^{\ap+k}E_{p-1-k}(-a)
  \mod p.$$\endpro
  Proof. For positive integers $m$ and $n$ it is well known ([MOS]) that
$$\sum_{r=0}^{m-1}(-1)^rr^n=\f{E_n(0)-(-1)^mE_n(m)}2.$$
Thus,
$$\sum_{r=1}^{\ap}\f{(-1)^r}{r^k}\e \sum_{r=0}^{\ap} (-1)^rr^{p-1-k}
=\f{E_{p-1-k}(0)-(-1)^{\ap+1}E_{p-1-k}(\ap+1)}2\mod p.$$ From [MOS]
and [S6, (2.2)-(2.3)] we know that
$$E_n(0)=\f{2(1-2^{n+1})B_{n+1}}{n+1}\qtq{and}
E_n(1-x)=(-1)^nE_n(x).\tag 2.2$$ Hence,
$$\sum_{r=1}^{\ap}\f{(-1)^r}{r^k}
\e -\f{(2^{p-k}-1)B_{p-k}}{p-k}+\f 12(-1)^{\ap+k}E_{p-1-k}(-\ap)\mod
p.$$ Set $a=\ap+pt$. It is well known ([MOS]) that
 $E_n(x+y)=\sum_{s=0}^n\b
ns x^sE_{n-s}(y)$. Thus,
$$\align E_{p-1-k}(-\ap)&=E_{p-1-k}(pt-a)=\sum_{s=0}^{p-1-k}
\b{p-1-k}s(pt)^sE_{p-1-k-s}(-a)\\&\e E_{p-1-k}(-a)\mod p.\endalign$$
We are done.
 \pro{Theorem 2.1} Let $p>3$ be a prime and
$a\in\Bbb Z_p$ with $a\not\e 0\mod p$. Then
$$\align \sum_{k=0}^{p-1}\b ak\b{-1-a}k&\e (-1)^{\ap}
+(a-\ap)(p+a-\ap)E_{p-3}(-a)
\\&\e (-1)^{\ap}
+(a-\ap)(p+a-\ap)\Big(\f 2{a^2}-E_{p-3}(a)\Big)
\mod{p^3}.\endalign$$
\endpro
Proof. Set $S_{p-1}(x)=\sum_{k=0}^{p-1}\b xk\b{-1-x}k$. Then
$$S_{p-1}(a)-(-1)^{\ap}S_{p-1}(a-\ap)
=\sum_{k=0}^{\ap-1}(-1)^k(S_{p-1}(a-k)+S_{p-1}(a-k-1)).$$ Suppose
$a=\ap+pt$. Then $t\in\Bbb Z_p$ and $a-k=\ap-k+pt$. For
$k=0,1,\ldots,\ap-1$ taking $m=\ap-k$ and $b=-1$ in [S8, (4.3)] we
see that
$$(a-k)(S_{p-1}(a-k)+S_{p-1}(a-k-1))
\e 2pt\cdot \f{-p-pt}{-(\ap-k)}=2p^2\f{t(t+1)}{\ap-k}\mod{p^3}$$ and
so
$$S_{p-1}(a-k)+S_{p-1}(a-k-1)\e 2p^2t(t+1)\cdot \f
1{(\ap-k)^2}\mod{p^3}.$$ Therefore,
$$\align S_{p-1}(a)-(-1)^{\ap}S_{p-1}(pt)
&=\sum_{k=0}^{\ap-1}(-1)^k(S_{p-1}(a-k)+S_{p-1}(a-k-1)) \\&\e
\sum_{k=0}^{\ap-1}(-1)^k\cdot 2p^2t(t+1)\cdot \f 1{(\ap-k)^2}
\\&=(-1)^{\ap}2p^2t(t+1)\sum_{r=1}^{\ap}\f{(-1)^r}{r^2}\mod{p^3}.
\endalign$$
As $B_{2m+1}=0$ for $m\ge 1$, we see that $B_{p-2}=0$. Thus, by
Lemma 2.2 we have
$$\sum_{r=1}^{\ap}\f{(-1)^r}{r^2}\e \f 12(-1)^{\ap}E_{p-3}(-a)\mod p.$$ Now, from the above and Lemma 2.1 we
deduce that
$$\align S_{p-1}(a)&\e
(-1)^{\ap}S_{p-1}(pt)+(-1)^{\ap}2p^2t(t+1)\sum_{r=1}^{\ap}\f{(-1)^r}{r^2}
\\&\e (-1)^{\ap}+p^2t(t+1)E_{p-3}(-a)\mod{p^3}.
\endalign$$
It is well known that ([MOS]) $E_n(1-x)=(-1)^nE_n(x)$ and
$E_n(x)+E_n(x+1)=2x^n$. Thus,
$$E_{p-3}(-a)=E_{p-3}(1+a)=2a^{p-3}-E_{p-3}(a)\e \f
2{a^2}-E_{p-3}(a)\mod p.$$
 Recall that $t=(a-\ap)/p$. By the above, the theorem is proved.

\par\q
\par Taking $a=-\f 12$ in Theorem 2.1 and then applying (1.11) and
the fact $E_n=2^nE_n(\f 12)$ we obtain (1.5).
\par For $m=3,4,6$ it is clear that
$$-\f 1m-\langle -\f 1m\rangle_p=\cases -\f 1m-\f{p-1}m=-\f pm&\t{if
$p\e 1\mod m$,}
\\-\f 1m-\f{(m-1)p-1}m=-\f{(m-1)p}m&\t{if $p\e -1\mod m$}
\endcases\tag 2.3$$ and so
$$\Big(-\f 1m-\langle -\f 1m\rangle_p\Big)
\Big(p-\f 1m-\langle -\f 1m\rangle_p\Big)=-\f pm\cdot
\f{(m-1)p}m=-\f {m-1}{m^2}p^2.\tag 2.4$$

\pro{Corollary 2.1} Let $p>3$ be a prime. Then
$$\sum_{k=0}^{p-1}\f{\b{6k}{3k}\b{3k}k}{432^k} \e \Ls {-1}p-\f
{25}9p^2E_{p-3}\mod{p^3}.$$
\endpro
Proof. Taking $a=-\f 16$ in Theorem 2.1 and then applying (1.11) and
(2.4) we see that
$$\align \sum_{k=0}^{p-1}\f{\b{6k}{3k}\b{3k}k}{432^k}
&=\sum_{k=0}^{p-1} \b{-\f 16}k\b{-\f 56}k
\\&\e (-1)^{\langle -\f 16\rangle_p}+\Big(-\f 16-\langle -\f 16\rangle_p\Big)
\Big(p-\f 16-\langle -\f 16\rangle_p\Big)E_{p-3}\Ls 16
\\&\e \Ls{-1}p-\f 5{36}E_{p-3}\Ls 16\mod{p^3}.
\endalign$$
  By [S6, Theorem 2.1 and
Lemma 2.1], we have $6^{2n}E_{2n}\sls 16=\f{3^{2n}+1}2E_{2n}.$ Thus,
$E_{p-3}(\f 16)=\f 1{6^{p-3}}\cdot \f{3^{p-3}+1}2E_{p-3}\e 20E_{p-3}
\mod p$. Hence the result follows.
\par\q
In [S7] the author introduced the sequence $\{U_n\}$ given by
$$U_0=1,\ U_n=-2\sum_{k=1}^{[n/2]}\b n{2k}U_{n-2k}\ (n\ge 1)$$
or
$$\sum_{n=0}^{\infty}U_n\f{t^n}{n!}=\f 1{\t{e}^t+\t{e}^{-t}-1}
\ (|t|<\f{\pi}3).$$ Clearly $U_{2n-1}=0$. The first few values of
$U_{2n}$ are shown below:
$$\align &U_2=-2,\q U_4=22,\q U_6=-602,\q U_8=30742,\q U_{10}=-2523002,
\\&U_{12}=303692662,\q U_{14}=-50402079002, \q U_{16}=11030684333782.
\endalign$$
For any prime $p>3$, in [S7] the author proved that
$$\sum_{k=1}^{[2p/3]}\f{(-1)^{k-1}}k\e 3p\Ls p3U_{p-3}\mod {p^2}.$$
  \pro{Corollary 2.2} Let $p>3$ be a prime. Then
$$\sum_{k=0}^{p-1}\f{\b{2k}k\b{3k}k}{27^k} \e \Ls {-3}p-2p^2U_{p-3}\mod{p^3}.$$
\endpro
Proof. Taking $a=-\f 13$ in Theorem 2.1 and then applying (1.11) and
(2.4) we see that
$$\align \sum_{k=0}^{p-1}\f{\b{2k}k\b{3k}k}{27^k}
&=\sum_{k=0}^{p-1} \b{-\f 13}k\b{-\f 23}k
\\&\e (-1)^{\langle -\f 13\rangle_p}+\Big(-\f 13-\langle -\f 13\rangle_p\Big)
\Big(p-\f 13-\langle -\f 13\rangle_p\Big)E_{p-3}\Ls 13
\\&=\Ls{-3}p-\f 29E_{p-3}\Ls 13\mod{p^3}.
\endalign$$
 By [S7, Theorem 2.1], $U_{2n}=3^{2n}E_{2n}\sls 13$. Thus,
$U_{p-3}=3^{p-3}E_{p-3}\sls 13\e \f 19E_{p-3}(\f 13)\mod p$. Now
putting all the above together we obtain the result.
\par\q
\newline{\bf Remark 2.1} Let $p>3$ be a prime. By [S7, p.217],
$B_{p-2}(\f 13)\e 6U_{p-3}\mod p$. Thus, from Corollary 2.2 we
deduce (1.8). In [MT], Mattarei and Tauraso proved that
$$\sum_{k=0}^{p-1}\b{2k}k\e \Ls{-3}p-\f{p^2}3B_{p-2}\Ls 13\mod
{p^3}.$$ This together with Corollary 2.2 yields
$$\sum_{k=0}^{p-1}\f{\b{2k}k\b{3k}k}{27^k}
\e\sum_{k=0}^{p-1}\b{2k}k\e \Ls{-3}p-2p^2U_{p-3}\mod{p^3}.$$
\par In [S3] the author introduced the sequence $\{S_n\}$ given by
$$S_0=1\qtq{and} S_n=1-\sum_{k=0}^{n-1}\b nk2^{2n-2k-1}S_k\ (n\ge 1),
$$
and showed that $S_n=4^nE_n(\f 14)$.
 \pro{Corollary 2.3} Let $p>3$
be a prime. Then
$$\sum_{k=0}^{p-1}\f{\b{2k}k\b{4k}{2k}}{64^k}\e \Ls{-2}p-3p^2S_{p-3}
\mod {p^3}.$$
\endpro
Proof. Taking $a=-\f 14$ in Theorem 2.1 and then applying (1.11) and
(2.4) we see that
$$\align \sum_{k=0}^{p-1}\f{\b{2k}k\b{4k}{2k}}{64^k}
&=\sum_{k=0}^{p-1} \b{-\f 14}k\b{-\f 34}k
\\&\e (-1)^{\langle -\f 14\rangle_p}+\Big(-\f 14-\langle -\f 14\rangle_p\Big)
\Big(p-\f 14-\langle -\f 14\rangle_p\Big)E_{p-3}\Ls 14
\\&=\Ls{-2}p-\f 3{16}E_{p-3}\Ls 14\mod{p^3}.
\endalign$$
Since $S_{p-3}=4^{p-3}E_{p-3}\sls 14\e  \f 1{16}E_{p-3}(\f 14)\mod
p$, we obtain the result.

\pro{Corollary 2.4} Let $p>3$ be a prime and $a\in\Bbb Z_p$ with
$a\not\e 0\mod p$. Then
$$\sum_{k=0}^{p-1}\b ak\b{-1-a}k+\sum_{k=0}^{p-1}\b{-a}k\b{-1+a}k \e
(a-\ap)(p+a-\ap)\f 2{a^2}\mod{p^3}.$$
\endpro
Proof. As $\langle -a\rangle_p=p-\ap$, from Theorem 2.1 we derive
that
$$\align \sum_{k=0}^{p-1}\b{-a}k\b{-1+a}k
&\e (-1)^{\langle -a\rangle_p} +(-a-\langle -a\rangle_p)(p-a-\langle
-a\rangle_p)\Big(\f 2{a^2}-E_{p-3}(-a)\Big)
\\&=-(-1)^{\ap}
+(a-\ap)(p+a-\ap)\Big(\f 2{a^2}-E_{p-3}(-a)\Big)
\\&\e (a-\ap)(p+a-\ap)\f 2{a^2}-\sum_{k=0}^{p-1}\b ak\b{-1-a}k \mod{p^3}.
\endalign$$
This yields the result.

\pro{Lemma 2.3} For any nonnegative integer $n$ we have
$$\sum_{k=0}^n(k-a(a+1))\b ak\b{-1-a}k=-a(a+1)\b{a-1}n\b{-2-a}n.$$
\endpro
Proof. Observe that
$$\align
&-a(a+1)\Big\{\b{a-1}{n+1}\b{-2-a}{n+1}-\b{a-1}n\b{-2-a}n\Big\}
\\&=\b a{n+1}\b{-1-a}{n+1}((a-n-1)(-2-a-n)-(n+1)^2)
\\&=(n+1-a(a+1))\b a{n+1}\b{-1-a}{n+1}.\endalign$$
The result can be easily proved by induction on $n$.

\pro{Theorem 2.2} Let $p>3$ be a prime and $a\in\Bbb Z_p$ with
$a\not\e 0,-1\mod p$. Then
$$\sum_{k=0}^{p-1}k\b ak\b{-1-a}k
\e (-1)^{\ap}a(a+1)+p^2t(t+1)\big(a(a+1)E_{p-3}(-a)-1\big)
\mod{p^3},$$ where $t=(a-\ap)/p$.
\endpro
Proof. By [S8, Lemma 4.2],
$$\align &\b{a-1}{p-1}=\b{\ap+pt-1}{p-1}\e \f {pt}{\ap}
\mod{p^2},
\\&\b{-2-a}{p-1}=\b{p-1-\ap-p(t+1)-1}{p-1}
\e \f{p(-t-1)}{p-1-\ap}\e \f{p(t+1)}{\ap+1}\mod{p^2}.
\endalign$$
Thus,
$$\b{a-1}{p-1}\b{-2-a}{p-1}\e \f{t(t+1)}{\ap(\ap+1)}p^2
\e \f{t(t+1)}{a(a+1)}p^2\mod{p^3}.$$ Hence, using Lemma 2.3 we see
that
$$\align &\sum_{k=0}^{p-1}k\b ak\b{-1-a}k
-a(a+1)\sum_{k=0}^{p-1}\b ak\b{-1-a}k
\\&=-a(a+1)\b{a-1}{p-1}\b{-2-a}{p-1}
\e -p^2t(t+1)\mod{p^3}.\endalign$$ This together with Theorem 2.1
yields the result.

\pro{Corollary 2.5} Let $p>3$ be a prime. Then
$$\sum_{k=0}^{p-1}\f{k\b{6k}{3k}\b{3k}k}{432^k}
\e -\f 5{36}\Ls{-1}p+p^2\Big(\f
5{36}+\f{125}{324}E_{p-3}\Big)\mod{p^3}.$$
\endpro
Proof. Taking $a=-\f 16$ in Theorem 2.2 and then applying (1.11) we
see that
$$\sum_{k=0}^{p-1}\f{k\b{6k}{3k}\b{3k}k}{432^k}
\e -\f 5{36}\Ls{-1}p-\f 5{36}p^2\Big(-\f 5{36}E_{p-3}\Ls
16-1\Big)\mod {p^3}.$$ By the proof of Corollary 2.1, $E_{p-3}\sls
16\e 20E_{p-3}\mod p$. Thus the result follows.

 \pro{Corollary 2.6}
Let $p>3$ be a prime. Then
$$\sum_{k=0}^{p-1}\f{k\b{2k}k\b{3k}k}{27^k}
\e -\f 29\Ls{-3}p+p^2\Big(\f 29+\f 49U_{p-3}\Big)\mod{p^3}.$$
\endpro
Proof. Taking $a=-\f 13$ in Theorem 2.2 and then applying (1.11) and
the fact $E_{p-3}(\f 13)\e 9U_{p-3}\mod p$ we deduce the result.

\subheading {3. Congruences for $\sum_{k=0}^{p-1}\b ak\b{-1-a}k\f 1
{2k+1}\mod {p^3}$}
\par For any positive integer $n$ and variables $a$ and $b$ with $
b\not\in \{-1,-\f 12,\ldots,-\f 1n\}$ let
$$S_n(a,b)=\sum_{k=0}^n\b ak\b{-1-a}k\f 1{bk+1}.\tag 3.1$$
Then
$$\align &(ab+1)S_n(a,b)-(ab-1)S_n(a-1,b)
\\&=\sum_{k=0}^n\b ak\b{-1-a}k\f {ab+1}{bk+1}-
\sum_{k=0}^n\b {a-1}k\b{-a}k\f {ab-1}{bk+1}
\\&=\sum_{k=0}^n\b ak\b{-a}k\Big(\f{ab+1}{bk+1}\cdot\f{a+k}a-
\f{ab-1}{bk+1}\cdot \f{a-k}a\Big)
\\&=2\sum_{k=0}^n\b ak\b{-a}k.
\endalign$$
By [S8, (4.5)] or induction on $n$,
$$\sum_{k=0}^n\b ak\b{-a}k=\b{n+a}n\b{n-a}n=\b{a-1}n\b{-a-1}n.$$
Thus,
$$(ab+1)S_n(a,b)-(ab-1)S_n(a-1,b)=2\b{a-1}n\b{-a-1}n.\tag 3.2$$

 \pro{Lemma 3.1 ([S8, Lemma 4.2])} Let $p$ be
an odd prime, $m\in\{1,2,\ldots,p-1\}$ and $t\in\Bbb Z_p$. Then
$$\b{m+pt-1}{p-1}\e
 \f{pt}m-\f{p^2t^2}{m^2}+\f{p^2t}mH_m
\mod{p^3}.$$\endpro

\pro{Lemma 3.2} Let $p$ be an odd prime and $a\in\Bbb Z_p$ with
$a\not\e 0\mod p$. Then
 $$\b{a-1}{p-1}\b{-a-1}{p-1}
 \e p^2\f{t(t+1)}{\ap^2}+p^3t(t+1)
\Big(-\f{1+2t}{a^3}+2\f{H_{\ap}}{a^2}\Big) \mod {p^4},$$
 where $t=(a-\ap)/p$.
 \endpro
Proof. By Lemma 3.1,
$$\align \b{a-1}{p-1}=\b{\ap+pt-1}{p-1}
&\e \f{pt}{\ap}+p^2t\Big(-\f t{\ap^2}+\f
1{\ap}H_{\ap}\Big)
\\&\e \f{pt}{\ap}+p^2t\Big(-\f t{a^2}+\f
{H_{\ap}}a\Big) \mod{p^3}.\endalign$$ From [S8, p.312] we know that
$H_{p-1-\ap}\e H_{\ap}\mod p$. Thus, from Lemma 3.1 we deduce that

$$\align \b{-a-1}{p-1}&=\b{p-\ap-p(t+1)-1}{p-1}
\\&\e \f{p(-t-1)}{p-\ap}+p^2(-t-1)\Big(-\f{-t-1}{(p-\ap)^2}
+\f{H_{p-\ap}}{p-\ap}\Big)
\\&\e
\f{p(t+1)(\ap+p)}{\ap^2}-p^2(t+1)
\Big(\f{t+1}{\ap^2}-\f{H_{p-\ap}}{\ap}\Big)
\\&\e \f{p(t+1)}{\ap}+p^2(t+1)\Big\{\f 1{\ap^2}-\f{t+1}{\ap^2}
+\f{-\f 1{\ap}+H_{p-1-\ap}}{\ap}\Big\}
\\&\e \f{p(t+1)}{\ap}+p^2(t+1)\Big\{-\f
{1+t}{\ap^2}+\f{H_{\ap}}{\ap}\Big\}
\\&\e \f{p(t+1)}{\ap}+p^2(t+1)\Big(-\f
{1+t}{a^2}+\f{H_{\ap}}{a}\Big) \mod{p^3}.\endalign$$ Hence,
$$\align &\b{a-1}{p-1}\b{-a-1}{p-1}
\\&\e \Big(\f{pt}{\ap}+p^2t\Big(-\f t{a^2}+\f
{H_{\ap}}a\Big)\Big)\Big( \f{p(t+1)}{\ap}+p^2(t+1)\Big(-\f
{1+t}{a^2}+\f{H_{\ap}}{a}\Big)\Big)
\\&\e p^2\f{t(t+1)}{\ap^2}+p^3t(t+1)
\Big(-\f{1+2t}{a^3}+2\f{H_{\ap}}{a^2}\Big)\mod{p^4}.
\endalign$$
This proves the lemma.
\par For any positive integer $n$ and variable $a$ let
$$T_n(a)=(2a+1)S_n(a,2)
=\sum_{k=0}^n\b ak\b{-1-a}k\f{2a+1}{2k+1}.\tag 3.3$$
 \pro{Lemma 3.3} Let $p>3$ be a prime and $t\in\Bbb Z_p$. Then
$T_{p-1}(pt)\e 1+2t\mod{p^3}.$
\endpro
Proof. Clearly
$$\align
T_{p-1}(pt)&=\sum_{k=0}^{p-1}\b{pt}k(-1)^k\b{pt+k}k\f{2pt+1}{2k+1}
\\&=2pt+1+\sum_{k=1}^{p-1}\f{(-1)^kpt(pt+k)(p^2t^2-(k-1)^2)\cdots (p^2t^2-1^2)}
{k!^2}\cdot \f{2pt+1}{2k+1}
\\&\e 2pt+1+\sum\Sb
k=1\\k\not=\f{p-1}2\endSb^{p-1}pt(pt+k)\f{(-1)^k(-1^2)(-2^2)\cdots
(-(k-1)^2)}{k!^2}\cdot\f{2pt+1}{2k+1}
\\&\qq+(-1)^{\f{p-1}2}\f{(2pt+1)t}{pt-\f{p-1}2}\cdot \f{(p^2t^2-\sls{p-1}2^2)\cdots
(p^2t^2-1^2)}{(\f{p-1}2!)^2}
\\&\e 2pt+ 1-pt(2pt+1)\sum\Sb
k=1\\k\not=\f{p-1}2\endSb^{p-1}\f{pt+k}{k^2(2k+1)}
+\f{2t(2pt+1)}{2pt+1-p} \Big(1-p^2t^2\sum_{k=1}^{\f{p-1}2}\f
1{k^2}\Big) \mod{p^3}.
\endalign$$
As $\f 1{k^2(2k+1)}=\f 1{k^2}-\f 2k+\f 4{2k+1}$, using (2.1) we see
that
 $$\align &\sum\Sb
k=1\\k\not=\f{p-1}2\endSb^{p-1}\f{pt+k}{k^2(2k+1)} =\sum\Sb
k=1\\k\not=\f{p-1}2\endSb^{p-1}(pt+k)\Big(\f 1{k^2}-\f 2k +\f
4{2k+1}\Big)
\\&=\sum_{k=1}^{p-1}(pt+k)\Big(\f 1{k^2}-\f 2k\Big)
-\Big(pt+\f{p-1}2\Big)\Big(\f 1{(\f{p-1}2)^2}-\f 2{\f{p-1}2}\Big)
+2\sum\Sb k=1\\k\not=\f{p-1}2\endSb^{p-1}\f{2pt-1+2k+1}{2k+1}
\\&\e pt\sum_{k=1}^{p-1}\Big(\f{1}{k^2}-2\f{1}k\Big)
+\sum_{k=1}^{p-1}\Big(\f{1}k-2\Big)-pt \Big(\f 1{\f 14}-\f 2{-\f
12}\Big)
\\&\qq-\Big(\f 1{\f{p-1}2}-2\Big)
+2\sum\Sb k=1\\k\not=\f{p-1}2\endSb^{p-1} 1 +2(2pt-1) \sum\Sb
k=1\\k\not=\f{p-1}2\endSb^{p-1}\f 1{2k+1}
\\&\e -8pt+2(p+1)+2(2pt-1) \sum\Sb
k=1\\k\not=\f{p-1}2\endSb^{p-1}\f 1{2k+1} \mod{p^2}.
\endalign$$
Also,
$$\align \sum\Sb
k=1\\k\not=\f{p-1}2\endSb^{p-1}\f{1}{2k+1}
&=\sum_{k=1}^{\f{p-3}2}\Big(\f 1{2k+1}+\f 1{2(p-k)+1} \Big)+\f
1{2\cdot \f{p+1}2+1}
\\&=\sum_{k=1}^{\f{p-3}2}\Big(\f 1{2k+1}+\f
{2p+2k-1}{(2p)^2-(2k-1)^2}\Big)+\f 1{p+2}
\\&\e -2p\sum_{k=1}^{\f{p-3}2}\f 1{(2k-1)^2}
+\sum_{k=1}^{\f{p-3}2} \Big(\f 1{2k+1}-\f 1{2k-1}\Big)+\f 1{p+2}
\\&=-2p\sum_{k=1}^{\f{p-1}2}\f 1{(2k-1)^2}+2p\cdot
\f 1 {(p-2)^2}-1+\f 1{p-2}+\f 1{p+2}
\\&\e -2p\sum_{k=1}^{\f{p-1}2}\f 1{(2k-1)^2}-1\mod{p^2}.
\endalign$$
By [S3, Corollary 2.1],
$$\align\sum_{k=1}^{\f{p-1}2}\f 1{(2k-1)^2}
&=\sum\Sb x=1\\x\e 1\mod 2\endSb^{p-1}\f 1{x^2} \e \sum\Sb x=0\\x\e
1\mod 2\endSb^{p-1}x^{p-3}
\\&\e  \f{2^{p-3}}{p-2}\big(B_{p-2}(0)-B_{p-2}(0)\big)=0\mod p.
\endalign$$
Hence,
$$\sum\Sb
k=1\\k\not=\f{p-1}2\endSb^{p-1}\f 1{2k+1}\e -1\mod{p^2}.\tag 3.4$$
Therefore,
$$\sum\Sb
k=1\\k\not=\f{p-1}2\endSb^{p-1}\f{pt+k}{k^2(2k+1)} \e
-8pt+2(p+1)-2(2pt-1)=2p(1-6t)+4 \mod{p^2}.$$ By [S2, Corollary 5.2],
$\sum_{k=1}^{(p-1)/2}\f 1{k^2}\e 0\mod p$. Thus, from all the above
we deduce that
$$\align T_{p-1}(pt)
&\e 2pt+1-pt(2pt+1)\sum\Sb
k=1\\k\not=\f{p-1}2\endSb^{p-1}\f{pt+k}{k^2(2k+1)}
+\f{2t(2pt+1)}{2pt+1-p}
\\&\e 2pt+1-pt(2pt+1)(2p(1-6t)+4)+2t\Big(1+\f{p}{1+(2t-1)p}\Big)
\\&\e 2pt+1-pt(8pt+2p(1-6t)+4)+2t+2tp(1-(2t-1)p)
\\&=1+2t\mod{p^3}.
\endalign$$
This proves the lemma.

 \pro{Theorem 3.1} Let $p>3$ be a prime and
$a\in\Bbb Z_p$ with $a(2a+1)\not\e 0\mod p$. Then
$$\align &\sum_{k=0}^{p-1}\b ak\b{-1-a}k\f 1{2k+1}
\\&\e \f{1+2t}{1+2a}
+p^2\f{t(t+1)}{1+2a}B_{p-2}(-a) \e \f{1+2t}{1+2a}
+p^2\f{t(t+1)}{1+2a}\Big(\f
2{a^2}-B_{p-2}(a)\Big)\mod{p^3},\endalign$$
 where $t=(a-\ap)/p$.
\endpro
Proof. As $a\not\e \f{p-1}2\mod p$, we see that
$$\align \b a{\f{p-1}2}\b{-1-a}{\f{p-1}2}&=(-1)^{\f{p-1}2}
\b a{\f{p-1}2}\b {a+\f{p-1}2}{\f{p-1}2} \\&=(-1)^{\f{p-1}2}
\f{(a+\f{p-1}2)(a+\f{p-1}2-1)\cdots (a-\f{p-1}2+1)}{(\f{p-1}2!)^2}
\e 0\mod p.\endalign$$ Thus, $\b ak\b{-1-a}k\f 1{2k+1}\in\Bbb Z_p$
for $k=0,1,\ldots,p-1$. Let $T_n(a)$ be given by (3.3). By (3.2) and
Lemma 3.2 we have
$$\aligned &T_{p-1}(a)-T_{p-1}(a-1)\\&=2\b{a-1}{p-1}\b{-a-1}{p-1}
\e 2p^2\f{t(t+1)}{\ap^2}+2p^3t(t+1)
\Big(-\f{1+2t}{a^3}+2\f{H_{\ap}}{a^2}\Big)\mod{p^4}.\endaligned \tag
3.5$$ For $1\le k\le \ap$ we have $\langle a-k+1\rangle_p=\ap-k+1$
and so $a-k+1=\ap-k+1+pt= \langle a-k+1\rangle_p+pt$. Hence
$$\align T_{p-1}(a)-T_{p-1}(a-\ap)&=\sum_{k=1}^{\ap}(T_{p-1}(a-k+1)-T_{p-1}(a-k))
\\&\e\sum_{k=1}^{\ap}\f { 2t(t+1)p^2}{\langle
a-k+1\rangle_p^2}=2t(t+1)p^2\sum_{k=1}^{\ap}\f 1{(\ap-k+1)^2}
\\&=2t(t+1)p^2\sum_{r=1}^{\ap}\f 1{r^2}\e 2t(t+1)p^2\sum_{r=1}^{\ap}
r^{p-3} \mod{p^3}.\endalign$$ By [S2, Lemma 3.2],
$$\sum_{r=1}^{\ap}
r^{p-3}\e (-1)^{p-2}\f{B_{p-2}(-a)-B_{p-2}}{p-2}\e \f 12B_{p-2}(-a)
\mod p.\tag 3.6$$ Thus,
$$\align T_{p-1}(a)-T_{p-1}(pt)&=T_{p-1}(a)-T_{p-1}(a-\ap)
\\&\e 2t(t+1)p^2
\cdot \f 12B_{p-2}(-a)=p^2t(t+1)B_{p-2}(-a)\mod {p^3}.\endalign$$
This together with Lemma 3.3 yields $T_{p-1}(a)\e
1+2t+p^2t(t+1)B_{p-2}(-a)\mod{p^3}.$ From [MOS] we know that
$B_n(-a)=(-1)^n(B_n(a)+na^{n-1})$. Thus,
$$B_{p-2}(-a)=(-1)^{p-2}(B_{p-2}(a)+(p-2)a^{p-3})\e -B_{p-2}(a)+
\f 2{a^2}\mod p.$$ This completes the proof.
 \pro{Corollary 3.1} Let
$p>3$ be a prime and $a\in\Bbb Z_p$ with $a\not\e 0,\pm\f 12\mod p$.
Then
$$\sum_{k=0}^{p-1}\b ak\b{-1-a}k\f{2a+1}{2k+1}
+\sum_{k=0}^{p-1}\b{-a}k\b{-1+a}k\f{1-2a}{2k+1} \e
(a-\ap)(p+a-\ap)\f 2{a^2}\mod{p^3}.$$
\endpro
Proof. As $\langle -a\rangle_p=p-\ap$, from Theorem 3.1 we derive
that
$$\align &\sum_{k=0}^{p-1}\b{-a}k\b{-1+a}k\f{1-2a}{2k+1}
\\&\e 1+2\f{-a-\langle -a\rangle_p}p +(-a-\langle
-a\rangle_p)(p-a-\langle -a\rangle_p)\Big(\f 2{a^2}-B_{p-2}(-a)\Big)
\\&=-1-2\f{a-\ap}p+(a-\ap)(p+a-\ap)\Big(\f 2{a^2}-B_{p-2}(-a)\Big)
\\&\e (a-\ap)(p+a-\ap)\f 2{a^2}
-\sum_{k=0}^{p-1}\b ak\b{-1-a}k\f{2a+1}{2k+1} \mod{p^3}.
\endalign$$
This yields the result.

\pro{Theorem 3.2} Let $p>3$ be a prime. Then
$$\sum_{k=0}^{p-1}\f{\b{2k}k\b{4k}{2k}}{64^k(2k+1)}
\e (-1)^{\f{p-1}2}-3p^2E_{p-3}\mod{p^3}.$$
\endpro
Proof. Taking $a=-\f 14$ in Theorem 3.1 and then applying (1.11),
(2.3) and (2.4) we obtain
$$\align \sum_{k=0}^{p-1}\f{\b{2k}k\b{4k}{2k}}{64^k(2k+1)}
&\e \f{1+2(-\f{2-\sls{-1}p}4)}{1+2(-\f 14)} +p^2\f{-\f 14\cdot\f 34}
{1+2(-\f 14)}B_{p-2}\Ls 14\\&=(-1)^{\f{p-1}2}-\f 38p^2B_{p-2}\Ls 14
\mod{p^3}.\endalign$$ It is known (see for example [S4, Lemma 2.5])
that $E_{2n}=-4^{2n+1}\f{B_{2n+1}(\f 14)}{2n+1}.$ Thus,
$$E_{p-3}=-4^{p-2}\f{B_{p-2}(\f 14)}{p-2}\e \f {B_{p-2}(\f 14)}8
\mod p.$$ Now combining all the above we obtain the result.

\pro{Theorem 3.3} Let $p>3$ be a prime. Then
$$\sum_{k=0}^{p-1}\f{\b{2k}k\b{3k}{k}}{27^k(2k+1)}
\e \Ls p3-4p^2U_{p-3}\mod{p^3}.$$
\endpro
Proof. Taking $a=-\f 13$ in Theorem 3.1 and then applying (1.11),
(2.3) and (2.4) we obtain
$$\sum_{k=0}^{p-1}\f{\b{2k}k\b{3k}{k}}{27^k(2k+1)}
\e \Ls p3-\f 23p^2B_{p-2}\Ls 13 \mod{p^3}.\tag 3.7$$ By [S7, p.217],
$B_{p-2}\sls 13\e 6U_{p-3}\mod p$. Thus the result follows.
\pro{Corollary 3.2} Let $p>3$ be a prime. Then
$$\sum_{k=0}^{p-1}\f{\b{2k}k\b{3k}{k}(4k+1)}{27^k(2k+1)}\e
\Ls p3\mod {p^3}.$$
\endpro
Proof. As $2-\f 1{2k+1}=\f{4k+1}{2k+1}$, combining Corollary 2.2
with Theorem 3.3 we deduce the result.

\pro{Theorem 3.4} Let $p>3$ be a prime. Then
$$\sum_{k=0}^{p-1}\f{\b{6k}{3k}\b{3k}{k}}{432^k(2k+1)}
\e \Ls p3-\f{25}4p^2U_{p-3}\mod{p^3}.$$
\endpro
Proof. Taking $a=-\f 16$ in Theorem 3.1 and then applying (1.11),
(2.3) and (2.4) we obtain
$$\sum_{k=0}^{p-1}\f{\b{6k}{3k}\b{3k}{k}}{432^k(2k+1)}
\e \Ls p3-\f 5{24}p^2B_{p-2}\Ls 16 \mod{p^3}.$$ By [S7, p.216],
$B_{p-2}\sls 16\e 30U_{p-3}\mod p$. Thus the result follows.
\par\q
\newline{\bf Remark 3.1} Corollary 3.2, (3.7) and the congruence
$\sum_{k=0}^{p-1}\f{\b{6k}{3k}\b{3k}{k}}{432^k(2k+1)} \e \sls
p3\mod{p^2}$ were conjectured by Z.W. Sun in [Su2].

\pro{Theorem 3.5} Let $p>3$ be a prime and $a\in\Bbb Z_p$ with
$2a+1\not\e 0\mod p$. Then
$$\sum_{k=0}^{p-1}\b ak\b{-1-a}k\f{a+k+1}{2k+1}
\e \f{1+(-1)^{\ap}}2+t+p^2\f{t(t+1)}4B_{p-2}\Big(-\f a2\Big) \mod
{p^3},$$ where $t=(a-\ap)/p$.
\endpro
Proof. From [MOS] we know that $E_n(x)=\f
2{n+1}(B_{n+1}(x)-2^{n+1}B_{n+1}\sls x2).$ Thus,
$$E_{p-3}(-a)=\f 2{p-2}\Big(B_{p-2}(-a)-2^{p-2}
B_{p-2}\Big(-\f a2\Big)\Big)\e -B_{p-2}(-a)+\f 12B_{p-2}\Big(-\f
a2\Big)\mod p.$$ Now from the above and Theorems 2.1 and 3.1 we
deduce that
$$\align &\sum_{k=0}^{p-1}\b ak\b{-1-a}k\Big(1+\f{2a+1}{2k+1}\Big)
\\&\e (-1)^{\ap}+p^2t(t+1)E_{p-3}(-a)+1+2t+p^2t(t+1)B_{p-2}(-a)
\\&\e 1+(-1)^{\ap}+2t+p^2\f{t(t+1)}2B_{p-2}\Big(-\f a2\Big)
\mod{p^3}.\endalign$$ This yields the result.

\subheading{4. Congruences for $\sum_{k=0}^{p-1}\b ak\b{-1-a}k\f
1{2k-1}\mod{p^3}$}

\pro{Lemma 4.1} For any nonnegative integer $n$ we have
$$\sum_{k=0}^n\b ak\b{-1-a}k
\f{(2a(a+1)+1)k-a(a+1)}{4k^2-1}=\f{a(a+1)}{2n+1}\b{a-1}n\b{-2-a}n.$$
\endpro
Proof. It is easy to check that
$$\align &\f{a(a+1)}{2(n+1)+1}\b{a-1}{n+1}\b{-2-a}{n+1}
-\f{a(a+1)}{2n+1}\b{a-1}n\b{-2-a}n
\\&=\b a{n+1}\b {-1-a}{n+1}\Big\{\f{a(a+1)}{2n+3}\cdot\f{a-n-1}a\cdot
\f{-2-a-n}{-1-a}-\f{a(a+1)}{2n+1}\cdot\f{n+1}a\cdot\f{n+1}{-1-a}\Big\}
\\&=\b a{n+1}\b {-1-a}{n+1}\f{(2a(a+1)+1)(n+1)-a(a+1)}{4(n+1)^2-1}.\endalign$$
Thus the result can be easily proved by induction on $n$.

\pro{Lemma 4.2} Let $p>3$ be a prime and $a\in\Bbb Z_p$ with
$a\not\e 0,-1,\pm \f 12\mod p$. Then
$$\sum_{k=0}^{p-1}\b ak\b{-1-a}k
\f{(2a(a+1)+1)k-a(a+1)}{4k^2-1}\e -(a-\ap)(p+a-\ap)\mod{p^3}.$$
\endpro
Proof. Set $a=\ap+pt$. By Lemma 3.1,
$$\align &\b{a-1}{p-1}=\b{\ap+pt-1}{p-1}\e \f{pt}{\ap}\e \f{pt}a\mod{p^2},
\\&\b{-2-a}{p-1}=\b{p-1-\ap-p(t+1)-1}{p-1}
\e \f{-p(t+1)}{p-1-\ap}\e \f{p(t+1)}{a+1}\mod{p^2}.\endalign$$ Thus,
$$\f{a(a+1)}{2(p-1)+1}\b{a-1}{p-1}\b{-2-a}{p-1}
\e \f{a(a+1)}{2p-1}\cdot\f {pt}a\cdot\f{p(t+1)}{a+1}\e -p^2t(t+1)
\mod{p^3}.$$ Now taking $n=p-1$ in Lemma 4.1 and then applying the
above we obtain the result.

 \pro{Theorem 4.1} Let $p>3$ be a prime and
$a\in\Bbb Z_p$ with $a\not\e 0,-1,\pm \f 12\mod p$. Then
$$\align&\sum_{k=0}^{p-1}\b ak\b{-1-a}k\f 1{2k-1}\\&\e
 -(2a+1)(2t+1)-p^2t(t+1)(4+(2a+1)B_{p-2}(-a)) \mod{p^3},
\endalign$$ where $t=(a-\ap)/p$.
\endpro
Proof. Note that $\f
1{2k-1}=4\f{(2a(a+1)+1)k-a(a+1)}{4k^2-1}-\f{(2a+1)^2}{2k+1}.$
Combining Theorem 3.1 with Lemma 4.2 we deduce the result.

\pro{Corollary 4.1} Let $p>3$ be a prime. Then
$$\sum_{k=0}^{p-1}\f{\b{2k}k\b{4k}{2k}}{64^k(2k-1)}
\e -\f 14\Ls{-1}p+\f 34p^2(1+E_{p-3})\mod{p^3}.$$
\endpro
Proof. Taking $a=-\f 14$ in Theorem 4.1 and noting that $B_{p-2}(\f
14)\e 8E_{p-3}\mod p$ we deduce the result.

\pro{Corollary 4.2} Let $p>3$ be a prime. Then
$$\sum_{k=0}^{p-1}\f{\b{2k}k\b{3k}{k}}{27^k(2k-1)}
\e -\f 19\Ls{-3}p+\f 49p^2(2-U_{p-3})\mod{p^3}.$$
\endpro
Proof. Taking $a=-\f 13$ in Theorem 4.1 and noting that $B_{p-2}(\f
13)\e 6U_{p-3}\mod p$ we deduce the result.

\pro{Corollary 4.3} Let $p>3$ be a prime. Then
$$\sum_{k=0}^{p-1}\f{\b{6k}{3k}\b{3k}{k}}{432^k(2k-1)}
\e -\f 49\Ls{-3}p+\f 59p^2(1+5U_{p-3})\mod{p^3}.$$
\endpro
Proof. Taking $a=-\f 16$ in Theorem 4.1 and noting that $B_{p-2}(\f
16)\e 30U_{p-3}\mod p$ we deduce the result.

\subheading {5. Congruences for $\sum_{k=1}^{p-1}\f 1k\b
ak\b{-1-a}k\mod {p^3}$}
\par For any positive integer $n$ and variable $a$ let
$$A_n(a)=\sum_{k=1}^n\b ak\b{-1-a}k\f 1k.$$
Then
$$\align A_n(a)-A_n(a-1)&=\sum_{k=1}^n\f 1k\Big\{\b
ak\b{-1-a}k-\b{a-1}k\b{-a}k\Big\}
\\&=\sum_{k=1}^n\f 1k\b ak\b{-a}k\Big(\f{a+k}a-\f{a-k}a\Big)
\\&=\f 2a\sum_{k=1}^n\b ak\b{-a}k.\endalign$$
By [S7, (4.5)] or induction on $n$, $$\sum_{k=0}^n\b ak\b{-a}k
=\b{n+a}n\b{n-a}n=\b{a-1}n\b{-a-1}n.$$ Thus,
$$A_n(a)-A_n(a-1)=\f 2a\b{a-1}n \b{-a-1}n-\f 2a.\tag 5.1$$
Hence, if $p>3$ is a prime and $a\in\Bbb Z_p$ with $a\not\e 0\mod
p$, by (5.1) and Lemma 3.2 we have
$$A_{p-1}(a)-A_{p-1}(a-1)\e \f{2t(t+1)}{a\ap^2}p^2-\f 2a\e
\f{2t(t+1)}{\ap^3}p^2-\f 2a\mod{p^3},\tag 5.2$$ where $t=(a-\ap)/p$.

\pro{Lemma 5.1} Let $p>3$ be a prime and $t\in\Bbb Z_p$. Then
$$\sum_{k=1}^{p-1}\b{pt}k\b{-1-pt}k\f 1k
\e -\f 23p^2tB_{p-3}\mod{p^3}.$$
\endpro
Proof. By the proof of Lemma 2.1, for $k\in\{1,2,\ldots,p-1\}$ we
have $\b{pt}k\b{-1-pt}k\e -\f{p^2t^2}{k^2}-\f{pt}k\mod{p^3}$. Thus,
$$\sum_{k=1}^{p-1}\b{pt}k\b{-1-pt}k\f 1k
\e -p^2t^2\sum_{k=1}^{p-1}\f 1{k^3}-pt\sum_{k=1}^{p-1}\f 1{k^2}
\mod{p^3}.$$ By [L] or [S2, Corollary 5.1], $\sum_{k=1}^{p-1}\f
1{k^3}\e 0\mod p$ and $\sum_{k=1}^{p-1}\f 1{k^2} \e \f
23pB_{p-3}\mod{p^2}.$ Thus the result follows.

\pro{Lemma 5.2} Let $p>3$ be a prime, $a\in\Bbb Z_p$, $a\not\e 0\mod
p$ and $t=(a-\ap)/p$. Then
$$\sum_{k=1}^{p-1}\f {\b ak\b{-1-a}k}k\e -\f 23p^2tB_{p-3}
-2\sum_{r=1}^{\ap}\f 1r+2pt\sum_{r=1}^{\ap}\f 1{r^2}
+2p^2t\sum_{r=1}^{\ap}\f 1{r^3}\mod{p^3}.$$
\endpro
Proof. For $1\le k\le \ap$ we have $\langle a-k+1\rangle_p=\ap-k+1$
and so $a-k+1=\langle a-k+1\rangle_p+pt$. Using (5.2) we see that
$$\align A_{p-1}(a)-A_{p-1}(a-\ap)&=\sum_{k=1}^{\ap}
(A_{p-1}(a-k+1)-A_{p-1}(a-k))
\\&\e\sum_{k=1}^{\ap}\Big(\f { 2t(t+1)p^2}{\langle
a-k+1\rangle_p^3}-\f 2{a-k+1}\Big) \\&=2t(t+1)p^2\sum_{k=1}^{\ap}\f
1{(\ap-k+1)^3}-2\sum_{k=1}^{\ap}\f 1{\ap-k+1+pt}
\\&=2t(t+1)p^2\sum_{r=1}^{\ap}\f 1{r^3}- 2
\sum_{r=1}^{\ap}\f 1{r+pt} \mod{p^3}.\endalign$$ Note that
$$\align \sum_{r=1}^{\ap}\f 1{r+pt}&
=\sum_{r=1}^{\ap}\f{r^2-ptr+p^2t^2}{r^3-(pt)^3} \e
\sum_{r=1}^{\ap}\f{r^2-ptr+p^2t^2}{r^3}
\\&=\sum_{r=1}^{\ap}\f 1r-pt\sum_{r=1}^{\ap}\f 1{r^2}
+p^2t^2\sum_{r=1}^{\ap}\f 1{r^3}\mod{p^3}.
\endalign$$
We then obtain
$$A_{p-1}(a)-A_{p-1}(a-\ap)\e
-2\sum_{r=1}^{\ap}\f 1r+2pt\sum_{r=1}^{\ap}\f 1{r^2}
+2tp^2\sum_{r=1}^{\ap}\f 1{r^3}\mod{p^3}.$$ By Lemma 5.1,
$A_{p-1}(a-\ap)=A_{p-1}(pt)\e -\f 23p^2tB_{p-3}\mod{p^3}$. Thus the
result follows.
\pro{Theorem 5.1} Let $p>3$ be a prime, $a\in\Bbb
Z_p$, $a\not\e 0\mod p$ and $t=(a-\ap)/p$. Then
$$\sum_{k=1}^{p-1}\f {\b ak\b{-1-a}k}k\e
-\f 23p^2t(t+1)B_{p-3}(-a) -2\f{B_{p^2(p-1)}(-a)-B_{p^2(p-1)}}
{p^2(p-1)} \mod{p^3}.$$
\endpro
Proof. It is well known that (see [MOS])
$$\align &\sum_{r=1}^mr^k=\f{B_{k+1}(m+1)-B_{k+1}}{k+1},\tag 5.3
\\&B_n(x+y)=\sum_{k=0}^n\b nk y^kB_{n-k}(x),\
B_n(1-x)=(-1)^nB_n(x).\tag 5.4\endalign$$ Thus, using Euler's
theorem we see that

$$\align &pt\sum_{r=1}^{\ap}\f 1{r^2}-\sum_{r=1}^{\ap}\f 1r
\\&\e pt\sum_{r=1}^{\ap}r^{p^2(p-1)-2}-\sum_{r=1}^{\ap}r^{p^2(p-1)-1}
\\&=pt\f{B_{p^2(p-1)-1}(\ap+1)-B_{p^2(p-1)-1}}{p^2(p-1)-1}-
\f{B_{p^2(p-1)}(\ap+1)-B_{p^2(p-1)}}{p^2(p-1)}
\\&=-pt\f{B_{p^2(p-1)-1}(-\ap)}{p^2(p-1)-1}-\f{B_{p^2(p-1)}(-\ap)
-B_{p^2(p-1)}}{p^2(p-1)}
\\&=-pt\f{B_{p^2(p-1)-1}(pt-a)}{p^2(p-1)-1}
-\f{B_{p^2(p-1)}(pt-a)-B_{p^2(p-1)}(-a)+B_{p^2(p-1)}(-a)-B_{p^2(p-1)}}
{p^2(p-1)}
\\&=
\f {-pt}{p^2(p-1)-1}\sum_{k=0}^{p^2(p-1)-1}\b{p^2(p-1)-1}k(pt)^k
B_{p^2(p-1)-1-k}(-a) \\&\q-\f
1{p^2(p-1)}\sum_{k=1}^{p^2(p-1)}\b{p^2(p-1)}k(pt)^kB_{p^2(p-1)-k}(-a)
-\f{B_{p^2(p-1)}(-a)-B_{p^2(p-1)}} {p^2(p-1)} \mod{p^3}.
\endalign$$
By [S1, Lemma 2.3], $B_m(-a)\in\Bbb Z_p$ for $m\not\e 0\mod{p-1}$
and $pB_m(-a)\in\Bbb Z_p$ for $m\e 0\mod{p-1}$. Thus,
$$\align &pt\sum_{r=1}^{\ap}\f 1{r^2}-\sum_{r=1}^{\ap}\f 1r
+\f{B_{p^2(p-1)}(-a)-B_{p^2(p-1)}} {p^2(p-1)}
\\&\e
pt\big(B_{p^2(p-1)-1}(-a)+(p^2(p-1)-1)ptB_{p^2(p-1)-2}(-a)\big)
\\&\qq-ptB_{p^2(p-1)-1}(-a)-\f{p^2(p-1)-1}2(pt)^2 B_{p^2(p-1)-2}(-a)
\\&\e -p^2t^2B_{p^2(p-1)-2}(-a)+\f 12p^2t^2B_{p^2(p-1)-2}(-a)
\\&=-\f 12p^2t^2B_{(p^2-1)(p-1)+p-3}(-a)\mod{p^3}.
\endalign$$
By [S2, Corollary 3.1],
$$B_{(p^2-1)(p-1)+p-3}(-a)\e ((p^2-1)(p-1)+p-3)
\f{B_{p-3}(-a)}{p-3}\e \f 23B_{p-3}(-a)\mod p.$$ Thus,
$$pt\sum_{r=1}^{\ap}\f 1{r^2}-\sum_{r=1}^{\ap}\f 1r
\e -\f{B_{p^2(p-1)}(-a)-B_{p^2(p-1)}} {p^2(p-1)}- \f
13p^2t^2B_{p-3}(-a)\mod{p^3}. \tag 5.5$$ By [S2, Lemma 3.2],
$$\sum_{r=1}^{\ap} \f 1{r^3}\e \sum_{r=1}^{\ap} r^{p-4}\e \f
{B_{p-3}(-a)-B_{p-3}}{p-3}\e -\f 13\big(B_{p-3}(-a)-B_{p-3}\big)
\mod p.$$ Thus, from Lemma 5.2 and (5.5) we derive that
$$\align \sum_{k=1}^{p-1}\f {\b ak\b{-1-a}k}k&\e-\f
23p^2tB_{p-3} -2\sum_{r=1}^{\ap}\f 1r+2pt\sum_{r=1}^{\ap}\f 1{r^2}
+2tp^2\sum_{r=1}^{\ap}\f 1{r^3}
\\&\e -\f
23p^2tB_{p-3}-2\f{B_{p^2(p-1)}(-a)-B_{p^2(p-1)}} {p^2(p-1)}
\\&\qq-\f 23p^2t^2B_{p-3}(-a)
-\f 23 p^2t\big(B_{p-3}(-a)-B_{p-3}\big)
\\&=-\f 23p^2t(t+1)B_{p-3}(-a)
-2\f{B_{p^2(p-1)}(-a)-B_{p^2(p-1)}} {p^2(p-1)} \mod{p^3}.\endalign$$
 This completes the proof.
 \pro{Lemma 5.3 ([MOS])} For any positive integer $n$ we have
 $$\align &B_{2n}\Ls 12=(2^{1-2n}-1)B_{2n},\ B_{2n}\Ls 13=\f{3-3^{2n}}{2\cdot 3^{2n}}B_{2n},\q
 \\&B_{2n}\Ls 14=\f{2-2^{2n}}{4^{2n}}B_{2n},\
 B_{2n}\Ls 16=\f{(2-2^{2n})(3-3^{2n})}{2\cdot 6^{2n}}B_{2n}.
 \endalign$$\endpro
 \par For an odd prime $p$ and $a\in\Bbb Z_p$ with $a\not\e 0\mod p$
 let $q_p(a)$ be the Fermat quotient given by
 $q_p(a)=(a^{p-1}-1)/p$. By Fermat's little theorem, $q_p(a)\in\Bbb
 Z_p$.

\pro{Theorem 5.2} Let $p>3$ be a prime. Then
$$\align &\sum_{k=1}^{p-1}\f{\b{2k}k\b{3k}k}{27^kk}
\e 3q_p(3)-\f
32pq_p(3)^2+p^2q_p(3)^3+\f{52}{27}p^2B_{p-3}\mod{p^3},\tag i
\\&\sum_{k=1}^{p-1}\f{\b{2k}k\b{4k}{2k}}{64^kk}
\e 6q_p(2)-3pq_p(2)^2+2p^2q_p(2)^3+\f 72p^2B_{p-3} \mod{p^3},\tag ii
\\&\sum_{k=1}^{p-1}\f{\b{6k}{3k}\b{3k}k}{432^kk} \e
4q_p(2)+3q_p(3)-p\Big(2q_p(2)^2 +\f 32q_p(3)^2\Big)\tag
iii\\&\qq\qq\qq+p^2\Big(\f
43q_p(2)^3+q_p(3)^3\Big)+\f{455}{54}p^2B_{p-3}
 \mod{p^3}.
\endalign$$
\endpro
Proof. By Lemma 5.3,
$$\align &B_{p-3}\Ls 13=\f{3-3^{p-3}}{2\cdot 3^{p-3}}B_{p-3}
\e 13B_{p-3}\mod p,
\\&B_{p-3}\Ls 14=\f{2-2^{p-3}}{4^{p-3}}B_{p-3}
\e 28B_{p-3}\mod p,
\\&B_{p-3}\Ls 16=\f{(2-2^{p-3})(3-3^{p-3})}{2\cdot 6^{p-3}}B_{p-3}
\e 91B_{p-3}\mod p.\endalign$$ By [S4, p.287],
$$\align &\f{B_{p^2(p-1)}-B_{p^2(p-1)}(\f 13)}{p^2(p-1)}
\e \f 32\big(q_p(3)-\f 12pq_p(3)^2+\f 13p^2q_p(3)^3\big)\mod
{p^3},\tag 5.6
\\&\f{B_{p^2(p-1)}-B_{p^2(p-1)}(\f 14)}{p^2(p-1)}
\e 3\big(q_p(2)-\f 12pq_p(2)^2+\f 13p^2q_p(2)^3\big)\mod {p^3},\tag
5.7
\\&\f{B_{p^2(p-1)}-B_{p^2(p-1)}(\f 16)}{p^2(p-1)}
 \e 2\Big(q_p(2)-\f 12pq_p(2)^2+\f 13p^2q_p(2)^3\Big)\tag 5.8
\\&\qq\qq\qq\qq\qq\q+\f 32\Big(q_p(3)-\f 12pq_p(3)^2+\f 13p^2q_p(3)^3\Big)\mod {p^3}.
\endalign$$
Now taking $a=-\f 13,-\f 14,-\f 16$ in Theorem 5.1 and then applying
(1.11) and the above we deduce the result.
\par\q
\newline{\bf Remark 5.1} Let $p>3$ be a prime. In [T3] Tauraso gave a congruence for
$\sum_{k=1}^{p-1}\f 1k\b ak\b{-1-a}k$ $\mod{p^2}$, and showed that
$\sum_{k=1}^{p-1}\f{\b{2k}k^2}{16^kk}\e -2H_{\f{p-1}2}\mod{p^3},$
which can be deduced from Theorem 5.1 (with $a=-\f 12$) and the
congruence ([S2, Theorem 5.2(c)])
$$H_{\f{p-1}2}\e -2q_p(2)+pq_p(2)^2-\f 23p^2q_p(2)^3-\f
7{12}p^2B_{p-3}\mod{p^3}.$$
 Theorem 5.2(ii) is equivalent to Z.W. Sun's
conjecture (1.10).
 \subheading {6. Congruences for
$\sum_{k=1}^{p-1}\f {(-1)^k}k\b ak\mod {p^2}$}
\par For given positive integer $n$ and variables $a$
 and $b$ define
$$f_n(a,b)=\sum_{k=1}^n\f{\b ak}{k\b bk}.$$
Then
$$f_n(a,b)-f_n(a-1,b)=\sum_{k=1}^n\f{\b ak-\b{a-1}k}{k\b bk}
=\sum_{k=1}^n\f{\f 1k\b{a-1}{k-1}}{\b bk}=\f 1a\sum_{k=1}^n\f{\b
ak}{\b bk}.$$ By Lerch's theorem ([B, p.86]) or induction on $n$,
$$\sum_{k=0}^n\f{\b
ak}{\b bk}=\f{b+1}{b+1-a}\Big\{1-\f{\b a{n+1}}{\b
{b+1}{n+1}}\Big\}.\tag 6.1$$ Thus,
$$\aligned f_n(a,b)-f_n(a-1,b)&=\f 1 a\Big\{\f{b+1}{b+1-a}-1-\f{b+1}{b+1-a}
\cdot\f{\b a{n+1}}{\b{b+1}{n+1}}\Big\}\\&=\f
1{b+1-a}-\f1{b+1-a}\cdot\f{\b{a-1}n}{\b bn}.\endaligned\tag 6.2$$
\pro{Lemma 6.1}Let $p>3$ be a prime, $n\in\{1,2,3,\ldots\}$,
$a,b\in\Bbb Z_p, 1\leq \ap\leq n\leq \langle b\rangle_p$ and
$a=\ap+pt$. Then
$$\aligned \sum_{k=1}^n\f{\b ak}{k\b bk}&\e pt\sum_{k=1}^n\f{(-1)^{k-1}}{k^2\b bk}+\sum_{r=1}^{\ap}
\f 1{b+1-r}+pt\sum_{r=1}^{\ap}\f 1{(b+1-r)^2}\\&\q-\f {pt}{\b
bn}\sum_{r=1}^{\ap}\f{(-1)^{n-r}}{(b+1-r)r\b
nr}\mod{p^2}.\endaligned$$

\endpro

Proof. By (6.2),
$$\aligned &f_n(a,b)-f_n(a-\ap,b)=\sum_{k=1}^{\ap}(f_n(a-k+1,b)-f_n(a-k,b))\\
&=\sum_{k=1}^{\ap}\f 1{b+1-(a-k+1)}-\sum_{k=1}^{\ap}\f
1{b+1-(a-k+1)}\cdot\f{\b{a-k+1-1}n}{\b bn}\\&=\sum_{k=1}^{\ap}\f
1{b+1-a+\ap-(\ap-k+1)}\\&\q-\f 1{\b bn}\sum_{k=1}^{\ap}\f
1{b+1-a+\ap-(\ap-k+1)}\b{\ap-k+1+a-\ap-1}n.\endaligned$$ Set
$a=\ap+pt$. Substituting $k$ with $\ap+1-r$ in the above we obtain
$$f_n(a,b)-f_n(pt,b)=\sum_{r=1}^{\ap}\f 1{b+1-r-pt}-\f1{\b bn}\sum_{r=1}^{\ap}
\f1{b+1-r-pt}\b{r+pt-1}n.\tag 6.3$$ For  $1\leq r\leq \ap$ we see
that
$$\aligned \b
{r-1+pt}n&=\f{(r-1+pt)(r-2+pt)\cdots(r-n+pt)}{n!}\\&=\f{(r-1+pt)(r-2+pt)\cdots(1+pt)
pt(pt-1)\cdots(pt-(n-r))}{n!}\\&\e\f{(r-1)!\cdot
pt\cdot(-1)^{n-r}\cdot(n-r)!}{n!}\\&=(-1)^{n-r}\f{pt}{r\b
nr}\mod{p^2}\endaligned$$ and so
$$\aligned f_n(a,b)-f_n(pt,b)&\e\sum_{r=1}^{\ap}\f 1{b+1-r-pt}-\f1{\b bn}
\sum_{r=1}^{\ap}\f{(-1)^{n-r}}{b+1-r-pt}\cdot\f{pt}{r\b
nr}\\&\e\sum_{r=1}^{\ap}\f{b+1-r+pt}{(b+1-r)^2}-\f{pt}{\b
bn}\sum_{r=1}^{\ap}\f{(-1)^{n-r}}{(b+1-r)r\b
nr}\mod{p^2}.\endaligned$$ On the other hand,$$
\sum_{k=1}^n\f{\b{pt}k}{k\b
bk}=\sum_{k=1}^n\f{pt\b{pt-1}{k-1}}{k^2\b bk}\e pt\sum_{k=1}^n\f{\b
{-1}{k-1}}{k^2\b bk}=pt\sum_{k=1}^n\f{(-1)^{k-1}}{k^2\b
bk}\mod{p^2}.$$ Thus, the result follows.

\pro{Theorem 6.1} Let $p>3$ be a prime, $a\in \Bbb Z_p$ and $a\not\e
0\mod p$. Then
$$\sum_{k=1}^{p-1}\f {(-1)^k}k\b ak\e
\f{B_{p^2(p-1)}-B_{p^2(p-1)}(-a)}{p^2(p-1)}
+\f{a-\ap}2B_{p-2}(-a)\mod{p^2}.$$
\endpro

Proof. Set $a=\ap+pt$. As $\b {-1}k=(-1)^k$, taking $b=-1$ and
$n=p-1$ in Lemma 6.1 and then applying the well known fact
$\sum_{k=1}^{p-1}\f 1{k^2}\e 0\mod p$, (5.5) and (3.6) we see that
$$\aligned \sum_{k=1}^{p-1}\f{(-1)^k}k\b
ak&=\sum_{k=1}^{p-1}\f{\b ak}{k\b{-1}k}\e-pt\sum_{k=1}^{p-1}\f
1{k^2}-\sum_{r=1}^{\ap}\f 1r+pt\sum_{r=1}^{\ap}\f
1{r^2}+pt\sum_{r=1}^{\ap}\f
1{r^2}\\&\e\f{B_{p^2(p-1)}-B_{p^2(p-1)}(-a)}{p^2(p-1)}
+\f{pt}2B_{p-2}(-a)\mod{p^2}.\endaligned$$
This proves the theorem.

\par\q
\newline{\bf Remark 6.1} In [T2] Tauraso showed that for any prime
$p>3$,
$$\sum_{k=1}^{p-1}\f{(-1)^k}k\b{-1/2}k\e -H_{\f{p-1}2}\mod{p^3}.$$
In [L] Lehmer proved that $H_{\f{p-1}2}\e
-2q_p(2)+pq_p(2)^2\mod{p^2}$. Thus,
$$\sum_{k=1}^{p-1}\f{(-1)^k}k\b{-1/2}k\e 2q_p(2)-pq_p(2)^2\mod{p^2}.$$
This can be deduced from Theorem 6.1 (with $a=-\f 12$) and Lemma
5.3. \pro{Theorem 6.2} Let $p>3$ be a prime. Then
$$\sum_{k=1}^{p-1}\f{(-1)^k}k\b{-1/4}k\e3q_p(2)-\f{3}2pq_p(2)^2-
\Big(2-\Ls{-1}p\Big) pE_{p-3}\mod{p^2}.$$\endpro
 Proof. Taking $a=-\f14$ in Theorem 6.1 we see
that
$$\sum_{k=1}^{p-1}\f{(-1)^k}k\b{-1/4}k\e\f{B_{p^2(p-1)}-B_{p^2(p-1)}(\f14)}{
p^2(p-1)}+\f{-\f 14-\langle -\f 14\rangle_p}2B_{p-2}\Ls
14\mod{p^2}.$$ By (5.7), $\f{B_{p^2(p-1)}-B_{p^2(p-1)}(\f14)}{
p^2(p-1)}\e 3q_p(2)-\f 32pq_p(2)^2\mod{p^2}.$ By the proof of
Theorem 3.2, $B_{p-2}\sls 14\e8E_{p-3}\mod p.$ Now, from the above
and (2.3) we deduce the result.

\pro{Theorem 6.3} Let $p>3$ be a prime. Then
$$\sum_{k=1}^{p-1}\f{(-1)^k}k\b{-1/3}k\e \f 32q_p(3)-\f 34pq_p(3)^2-
\f{3-\sls{-3}p}2pU_{p-3}\mod{p^2}.$$\endpro

Proof. Taking $a=-\f 13$ in Theorem 6.1 we see that
$$\sum_{k=1}^{p-1}\f{(-1)^k}k\b{-1/3}k\e
\f{B_{p^2(p-1)}-B_{p^2(p-1)}(\f 13)}{ p^2(p-1)}+\f{-\f 13-\langle
-\f 13\rangle_p}2B_{p-2}\Ls 13\mod{p^2}.$$ By (5.6),
$\f{B_{p^2(p-1)}-B_{p^2(p-1)}(\f 13)}{ p^2(p-1)}\e \f 32q_p(3)-\f
34pq_p(3)^2\mod{p^2}.$ By [S7, p.217], $B_{p-2}\sls 13\e
6U_{p-3}\mod p.$ Now, from the above and (2.3) we deduce the result.

\pro{Theorem 6.4} Let $p>3$ be a prime. Then
$$\align &\sum_{k=1}^{p-1}\f{(-1)^k}k\b{-1/6}k\\&\e 2q_p(2)+\f 32q_p(3)
-pq_p(2)^2-\f 34pq_p(3)^2-\f
52\Big(3-2\Ls{-3}p\Big)pU_{p-3}\mod{p^2}.\endalign$$\endpro

Proof. Taking $a=-\f 16$ in Theorem 6.1 and then applying (5.8),
(2.3) and the fact $B_{p-2}\sls 16\e 30U_{p-3}\mod p$ we deduce the
result.

 \subheading {7. Congruences for
$\sum_{k=0}^{p-1}\b ak(-2)^k\mod {p^2}$}
\par For given positive integer $n$ and variables $a$
 and $x$ define
$$F_n(a,x)=\sum_{k=0}^n\b ak{x^k}.$$
Then
$$\align  F_n(a,x)-(x+1)F_n(a-1,x)&=
\sum_{k=0}^n\Big(\b ak-\b{a-1}k\Big)x^k-\sum_{k=0}^n\b{a-1}kx^{k+1}
\\&=\sum_{k=1}^n\b{a-1}{k-1}x^k-\sum_{k=0}^n\b{a-1}kx^{k+1}.\endalign$$
Thus, $$F_n(a,x)-(x+1)F_n(a-1,x)=-\b{a-1}nx^{n+1}.\tag 7.1$$
Suppose
that $p>3$ is a prime, $a\in\Bbb Z_p$ and $a=\ap+pt$. Taking $n=p-1$
in (7.1) and then applying Lemma 3.1 we see that
$$\align &F_{p-1}(a,x)-(x+1)F_{p-1}(a-1,x)
\\&=-\b{\ap+pt-1}{p-1}x^p
\e \Big(-\f{pt}{\ap}+\f{p^2t^2}{\ap^2}-\f{p^2t}{\ap}H_{\ap}\Big)x^p
\mod{p^3}.\endalign$$ For $1\le k\le \ap$ we have $\langle
a-k+1\rangle_p=\ap-k+1$ and so $a-k+1=\langle a-k+1\rangle_p+pt$.
Thus,
$$\align &F_{p-1}(a,x)-(x+1)^{\ap}F_{p-1}(a-\ap,x)
\\&=\sum_{k=1}^{\ap}(x+1)^{k-1}(F_{p-1}(a-k+1,x)-(x+1)F_{p-1}(a-k,x))
\\&\e \sum_{k=1}^{\ap}(x+1)^{k-1}x^p
\Big(-\f{pt}{\ap -k+1}
+\f{p^2t^2}{(\ap-k+1)^2}-\f{p^2t}{\ap-k+1}H_{\ap-k+1}\Big)
\\&=x^p\sum_{r=1}^{\ap}(x+1)^{\ap-r}\Big(-\f{pt}r
+\f{p^2t^2}{r^2}-\f{p^2t}rH_r\Big)\mod{p^3}
\endalign$$
and therefore
$$\aligned &F_{p-1}(a,x)-(x+1)^{\ap}F_{p-1}(pt,x)\\&\e x^p(x+1)^{\ap}\Big(-pt\sum_{r=1}^{\ap}\f 1{r(x+1)^r}
+p^2t^2\sum_{r=1}^{\ap}\f 1{r^2(x+1)^r}-p^2t\sum_{r=1}^{\ap}\f
{H_r}{r(x+1)^r}\Big)\mod{p^3}.
\endaligned\tag 7.2$$
\pro{Theorem 7.1} Let $p>3$ be a prime, $a\in\Bbb Z_p$ and $a\not\e
0\mod p$. Then
$$\sum_{k=0}^{p-1}\b ak(-2)^k\e (-1)^{\ap}-(a-\ap)E_{p-2}(-a)
\mod{p^2}.$$
\endpro
Proof. Set $t=(a-\ap)/p$. Taking $x=-2$ in (7.2) we see that
$$\sum_{k=0}^{p-1}\b ak(-2)^k-(-1)^{\ap}
\sum_{k=0}^{p-1}\b {pt}k(-2)^k \e
2^p(-1)^{\ap}pt\sum_{r=1}^{\ap}\f{(-1)^r}r\mod{p^2}.$$ By a result
of Glaisher (see [S4]), $\sum_{k=1}^{p-1}\f{2^k}k\e -2q_p(2)\mod p$.
Note that $\b{pt-1}{k-1}\e \b{-1}{k-1}=(-1)^{k-1} \mod p$ for $1\le
k\le p-1$. We then derive that
$$\align \sum_{k=0}^{p-1}\b {pt}k(-2)^k&=1+pt\sum_{k=1}^{p-1}
\b{pt-1}{k-1}\f{(-2)^k}k
\\&\e 1-pt\sum_{k=1}^{p-1}\f{2^k}k\e 1+2ptq_p(2)\mod{p^2}.
\endalign$$
It is well known that $pB_{p-1}\e p-1\mod p$. Thus, from Lemma 2.2
we deduce that
$$\align \sum_{r=1}^{\ap}\f{(-1)^r}r&\e -\f{q_p(2)pB_{p-1}}{p-1}+\f 12
(-1)^{\ap+1}E_{p-2}(-a)\\&\e -q_p(2)-\f 12(-1)^{\ap}E_{p-2}(-a)\mod
p.\endalign$$ Now combining all the above we deduce that
$$\align \sum_{k=0}^{p-1}\b ak(-2)^k&\e (-1)^{\ap}(
1+2ptq_p(2))+2(-1)^{\ap}pt\Big(-q_p(2)-\f
12(-1)^{\ap}E_{p-2}(-a)\Big)
\\&=(-1)^{\ap}-ptE_{p-2}(-a)\mod{p^2}.
\endalign$$
This proves the theorem.

\pro{Theorem 7.2} Let $p>3$ be a prime. Then
$$\sum_{k=0}^{p-1}\b {-1/3}k(-2)^k\e \Ls{-3}p+\f{3-\sls{-3}p}3
pq_p(2) \mod{p^2}.$$
\endpro
Proof. Taking $a=-\f 13$ in Theorem 7.1 and applying (2.3) we see
that
$$\align \sum_{k=0}^{p-1}\b {-1/3}k(-2)^k&\e (-1)^{\langle -\f 13\rangle_p}
-\Big(-\f 13-\langle -\f 13\rangle_p\Big) E_{p-2}\Ls 13
\\&=\Ls{-3}p-\f{\sls{-3}p-3}6pE_{p-2}\Ls 13\mod{p^2}.
\endalign$$
By [S7, Lemma 2.2], Lemma 5.3 and the fact $pB_{p-1}\e p-1\mod p$,
$$\align E_{p-2}\Ls 13&=\f 2{p-1}((-2)^{p-1}-1)B_{p-1}\Ls 13
=\f 2{p-1}\cdot pq_p(2)\cdot\f{3-3^{p-1}}{2\cdot 3^{p-1}}B_{p-1}
\\&\e 2q_p(2)\mod p.\endalign$$
Thus the result follows.
\par\q
\newline{\bf Remark 7.1} In [Su1], Z.W. Sun proved that for any
prime $p>3$,
$$\sum_{k=1}^{p-1}\f{\b{2k}k}{2^k}=\sum_{k=1}^{p-1}\b{-1/2}k(-2)^k
\e \Ls{-1}p-p^2E_{p-3}\mod{p^3}.$$

\enddocument

\Refs\widestnumber\key{BEW}
 \ref\key B\by H. Bateman\book
Higher Transcendental Functions\publ Vol.I, McGraw-Hill\publaddr New
York\yr 1953\endref

\ref\key L\by E. Lehmer\paper On congruences involving Bernoulli
numbers and the quotients of Fermat and Wilson \jour Ann. of
Math.\vol 39\yr 1938\pages 350-360
\endref
\ref\key MOS\by W. Magnus, F. Oberhettinger and R.P. Soni\book
Formulas and Theorems for the Special Functions of Mathematical
Physics $(3rd\ ed.)$\publ Springer-Verlag, New York\yr 1966\pages
25-32\endref
\ref\key MT\by S. Mattarei and R. Tauraso \paper
Congruences for central binomial sums and finite polylogarithms
\jour J. Number Theory\vol 133\yr 2013\pages 131-157
\endref


  \ref\key M1\by
E. Mortenson\paper A supercongruence conjecture of
Rodriguez-Villegas for a certain truncated hypergeometric
function\jour J. Number Theory \vol 99\yr 2003 \pages 139-147\endref

 \ref\key M2\by  E. Mortenson\paper Supercongruences
between truncated $\ _2F_1$ hypergeometric functions and their
Gaussian analogs\jour  Trans. Amer. Math. Soc. \vol 355\yr
2003\pages 987-1007\endref


 \ref\key RV\by F. Rodriguez-Villegas
\paper  Hypergeometric families of Calabi-Yau manifolds\jour in:
Noriko Yui, James D. Lewis (Eds.), Calabi-Yau Varieties and Mirror
Symmetry , Toronto, ON, 2001, in: Fields Inst. Commun., vol. 38,
Amer. Math. Soc., Providence, RI, 2003, pp.223-231\endref




\ref\key S1\by Z.H. Sun \paper Congruences for Bernoulli numbers and
Bernoulli polynomials\jour Discrete  Math. \vol 163\yr 1997\pages
153-163\endref

\ref\key S2\by Z.H. Sun \paper Congruences concerning Bernoulli
numbers and Bernoulli polynomials\jour Discrete Appl. Math. \vol
105\yr 2000\pages 193-223\endref

\ref\key S3\by Z.H. Sun\paper Congruences involving Bernoulli
polynomials\jour  Discrete Math. \vol 308\yr 2008\pages
71-112\endref

\ref\key S4\by Z.H. Sun\paper Congruences involving Bernoulli and
Euler numbers\jour J. Number Theory\vol 128\yr 2008\pages 280-312
\endref

\ref\key S5\by Z.H. Sun\paper Congruences concerning Legendre
polynomials \jour Proc. Amer. Math. Soc. \vol 139\yr 2011\pages
1915-1929\endref


 \ref\key S6\by Z. H. Sun\paper Congruences for sequences
similar to Euler numbers\jour J. Number Theory \vol 132\yr 2012
\pages 675-700\endref

\ref\key S7\by Z. H. Sun\paper Identities and congruences for a new
sequence\jour Int. J. Number Theory \vol 8\yr 2012\pages
207-225\endref



\ref\key S8\by  Z.H. Sun\paper Generalized Legendre polynomials and
related supercongruences\jour J. Number Theory\vol 143\yr 2014
\pages 293-319
\endref


\ref\key Su1\by Z.W. Sun\paper Super congruences and Euler numbers
\jour Sci. China Math. \vol 54\yr 2011\pages 2509-2535\endref

\ref\key Su2\by Z.W. Sun\paper p-adic congruences motivated by
series\jour J. Number Theory\vol 134\yr 2014\pages 181-196\endref


\ref \key T1\by R. Tauraso\paper An elementary proof of a
Rodriguez-Villegas supercongruence\jour arXiv:0911.4261, 2009\endref

\ref \key T2\by R. Tauraso\paper Congruences involving alternating
multiple harmonic sums\jour Electron. J. Combin. \vol 17\yr
2010\pages \#R16, 11 pp
\endref


\ref \key T3\by R. Tauraso\paper Supercongruences for a truncated
hypergeometric series\jour Integers \vol 12\yr 2012\pages \#A45, 12
pp
\endref



\endRefs
\bye